\documentclass[a4paper]{article}

\usepackage[latin1]{inputenc} %
\usepackage[T1]{fontenc} %
\usepackage{RR}
\usepackage{hyperref}

\usepackage{graphicx,amsmath}
\usepackage{subfigure}
\usepackage{epsfig}
\usepackage{verbatim}
\usepackage{algorithm}
\usepackage{algorithmic}
\usepackage{color}
\usepackage{multirow}
\usepackage{latexsym}

\newcommand{\eps}{\varepsilon}

\newcommand{\cA}{\mathcal{A}}

\newcommand{\cH}{\mathcal{H}}

\def\Heav{\mathop{\rm Heav}\nolimits}

\newcommand{\parder}[2]{\frac{\partial #1}{\partial #2}}

\newcommand{\beq}{\begin{equation}}
\newcommand{\eeq}{\end{equation}}

\newcommand{\eg}{{\it e.g. }}

\RRdate{June 2010}

\RRauthor{%
A.R. Akhmetzhanov\thanks{INRIA Sophia Antipolis-M\'editerran\'ee, France, akhmetzhanov@gmail.com}
\and
F. Grognard\thanks{INRIA Sophia Antipolis-M\'editerran\'ee, France, Frederic.Grognard@sophia.inria.fr}
\and
L. Mailleret\thanks{INRA Sophia Antipolis, France, Ludovic.Mailleret@sophia.inra.fr}
\and
P. Bernhard\thanks{INRIA Sophia Antipolis-M\'editerran\'ee, France, Pierre.Bernhard@sophia.inria.fr}
}

\authorhead{Akhmetzhanov \& Grognard \& Mailleret \& Bernhard}

\RRtitle{Se rallier ou filouter ?\\ Une analyse évolutionnaire dans un système consommateurs-ressources. }
\RRetitle{Join forces or cheat: evolutionary analysis of a consumer-resource system}
\titlehead{Join forces or cheat: evolutionary analysis of a consumer-resource system}

\RRresume{Dans ce travail nous étudions les processus évolutifs à l'oeuvre dans un système consommateurs-ressources ``saisonnier" dans lequel les consommateurs ont un compromis dynamique à faire entre le temps alloué à la recherche de nourriture et celui alloué à leur reproduction. Les individus sont actifs pendant des saisons de longueur fixe, séparées par des périodes d'hiver où seuls les immatures produits durant la saison survivent (oeufs, graines,...).  La taille de la génération suivante de matures est alors déterminée par ce nombre de survivants.

Dans ce rapport,  nous considérons d'abord une population de consommateurs qui maxi\-mi\-sent leur ``fitness'' commune (l'efficacité reproductive de la population entière), tous les consommateurs ayant ce même objectif; ils agissent en quelque sorte pour le bien commun. Nous supposons par la suite qu'une petite fraction de la population des consommateurs peut apparaître au début d'une saison et choisir une stratégie différente de celle de la population principale; nous appellerons ces individus déviants  ``mutants". Nous étudions en détails la stratégie mise en oeuvre par ces mutants et analysons pour finir leur capacité ou non à supplanter la population résidente.
}
\RRabstract{
In this work we study the process of mutant invasion on an example of a consumer-resource system with annual character of the behavior. Namely, individuals are active during seasons of fixed length separated by winter periods. All individuals die at the end of the season and the size of the next generation is determined by the number of offspring produced during the past season. The rate at which the consumers produce immature offspring depends on their internal energy which can be increased by feeding. The reproduction of the resource simply occurs at a constant rate.

At the beginning, we consider a population of consumers maximizing their common fitness, all consumers being individuals having the same goal function and acting for the common good. We suppose that a small fraction of the consumer population may appear at the beginning of one season and start to behave as mutants in the main population. We study how such invasion occurs.
}

\RRmotcle{Commande optimale, Jeux différentiels, Modèle consommateur-resource, Invasion}
\RRkeyword{Optimal Theory, Differential Games, Resource-Consumer Model, Mutant Invasion}

\RRprojet{Comore}  
\RRtheme{\THBio} 
\URSophia 

\begin{document}

\RRNo{7312}

\makeRR

\section{Introduction}
\label{Section_Introduction}
Biodiversity found on Earth consists of millions of biological species, thousands of different ecosystems. Among this variety, one can easily identify many examples of resource-consumer systems like prey-predator/parasitoid-host systems known in biology \cite{Murray1989} or birth-death systems known in cell biochemistry \cite{Assaf2008}. Usually, individuals involved in such kind of systems (bacteria, plants, insects, animals) have conflicting interests and models describing such interactions are based on principles of game theory \cite{Perrin2000,Houston2005,Auger2006,Hamelin2007}. Hence the investigation of these models is of interest both for game theoreticians and for biologists working in behavioral and evolutionary ecology.

One of the first questions that usually appears when first consulting evolutionary theory books is the following: could we say that individuals behave rationally or optimally throughout their life? The answer is most probably ``yes'' if we consider the evolution as a slow process tending to some equilibrium. Following Darwin theory and its main statement about the survival of the fittest, we can assume that evolution of populations leads to a situation where individuals maximize their fitness or try to protect themselves from invasion by others \cite{MaynardSmith1982,Vincent2005}. Such population can be referred to as \emph{residents} who use an \emph{optimal maximizing strategy} or an \emph{uninvadable strategy} respectively. The first type of strategy could be dynamically stable and lead to an asymptotically stable equilibrium, but it could also not be. Particularly, and this is a well-known fact in economics, a \emph{free-rider} may overcompete competitors cooperating with him by ``cheating'' and using a ``greedy'' strategy. In the sequel, populations which behave differently from the residents will be termed \emph{mutants}. On the other hand, if residents use an evolutionary stable strategy, this will not allow them to get the maximum possible value of the fitness but will help them avoid mutant invasion. This seems reasonable from a biological point of view but such strategy could be dynamically unreachable or could not lead to a stable equilibrium in a long-term perspective \cite{Brown2008}.

In this work we study the process of mutant invasion on an example of a consumer-resource system with annual character of the behavior as introduced by \cite{Akhmetzhanov2010}. Namely, individuals are active during seasons of fixed length $T$ separated by \emph{winter periods}. To give a representation of what such a system could encompass, the resource population could represent plants producing seeds all season long, and the consumer population insects having to trade-off between \emph{feeding} and \emph{laying eggs}. All individuals die at the end of the season and the size of the next generation is determined by the number of offspring (seeds or eggs) produced during the past season. The rate at which the consumers produce immature offspring (\emph{eggs}) depends on their internal energy which can be increased by feeding. The reproduction of the resource simply occurs at a constant rate.

In nature several patterns of life-history of the consumers can be singled out, but they almost always contains two main phases: \emph{growth phase} and \emph{reproduction phase}. Depending on initial conditions the transition between them could be strict when the consumers only feed at the beginning of their life and only reproduce at the end, or there could exist  an \emph{intermediate phase} between them where growth and reproduction occur simultaneously. Such types of behaviors are called \emph{determinate} and \emph{indeterminate growth pattern} respectively \cite{Perrin1993}.

Time-sharing between laying eggs and feeding for the consumers is described by the variable $u$: $u=1$ means feeding, $u=0$ on the other side means reproducing. The intermediate control $u\in(0,1)$ describes a situation where, for some part of the time, the individual is feeding and, for the other part of the time, it is reproducing.

At the beginning of the paper, we consider a population of consumers maximizing their common fitness, all consumers being individuals having the same goal function and {\it acting for the common good}. We will call them the residents in the following. We suppose that a small fraction of the consumer population may appear at the beginning of one season and start to behave as mutants in the main population. We study how such invasion will occur.

If there is a large number of residents, it makes sense to assume that residents fix their strategy \emph{a priori} and do not change it during the season. The mutant achieves better result than the resident using this fact and react in feedback form. Such a problem can be related to a hierarchical game of two players. For simplicity we study this problem in the case of a vanishingly small number of mutants. Such a situation corresponds to the definition of an evolutionarily stable strategy given by \cite{MaynardSmith1982}, when only a small fraction of mutants is taken into consideration. We also investigate the fate of such a mutation in the multi-seasonal framework proposed by \cite{Akhmetzhanov2010}. In particular we show that mutants not only can invade the resident consumers' population, but will also replace it in the system. Finally, we make some conclusions regarding the results presented in the paper.

\section{Main model}

\subsection{Previous work}
\label{s:optcoll}

\begin{figure}[t]
$$\includegraphics[scale=.45]{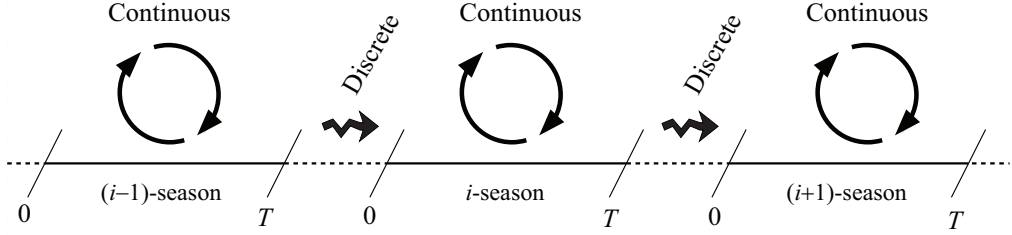}$$
\caption{Seasonal character of the behavior: populations evolve continuously during the season and have discrete rules for the transition from one season to another. Picture has been taken from \cite{Mailleret2009}}
\label{f:1}
\end{figure}

At the beginning consider a system of two populations: consumers and resource without any mutant. As it has been stated, all seasons have fixed length $T$ which does not change from one year to another (see Fig.~\ref{f:1}). The consumer population is determined by two state variables: the average energy of one individual $p$ and the number of consumers $c$ present in the system. For the description of the resource population a variable $n$ is introduced. It defines the size of the population. We suppose that both populations consist of two parts: \emph{mature} (insects/plants) and \emph{immature part} (eggs/seeds). During the season, mature individuals can invest in immatures by laying eggs. Between seasons (at winter periods) all matures die and immatures become matures for the next season.

We suppose that all consumers have arbitrarily small energy $p$ at the beginning of the season. The efficiency of the reproduction is assumed to be proportional to the value of $p$; it is thus intuitively understandable that consumers should feed on the resource at the beginning and reproduce at the end once they have gathered enough energy. The consumer has a trade-off between feeding ($u=1$) and laying eggs ($u=0$). The variable $u$ plays the role of the control.

The within season dynamics are thus defined as follows
\begin{equation}
\label{e:1}
\dot p = -ap+bnu,\quad \dot n = -cnu
\end{equation}
where we supposed that both populations do not suffer from mortality; $a$, $b$ and $c$ are some constants. After rescaling of time and state variables, the constants $a$ and $b$ can be eliminated and equations (\ref{e:1}) can be rewritten in a simplified form
\begin{equation}
\label{e:dynamics0}
\dot p = -p+nu,\quad \dot n = -cnu
\end{equation}
where $c$ is represent the number of predators present in the system.

The amount of offspring produced by individual during the season depends on the current size of the populations
\begin{equation}
\label{e:value0}
  J=\int_0^T \theta c(1-u(t))p(t)\>\mathrm dt,\quad J_n=\int_0^T \gamma n(t)\>\mathrm dt
\end{equation}
where consumers are maximizing the value $J$, the common fitness, $\theta$ and $\gamma$ are some constants. We see that this is an optimal control problem which can be solved using the dynamic programming \cite{Bellman1957} or Pontryagin maximum principle \cite{Pontryagin1962}. Moreover, the constants $c$, $\theta$ and $\gamma$ can be omitted to compute the solution of this problem without loss of generality.

One can also show that all the data of the formulated problem are homogeneous of degree one in state variables, which can be only positive numbers. This is a particular case of Noether's theorem in the calculus of variations about the problems whose data is invariant under a group of transformations \cite{Caratheodory1965}. Hence the dimension of phase space of the optimal control problem (\ref{e:dynamics0}-\ref{e:value0}) can be lowered by one unit by the introduction of a new variable $x=p/n$. In this case its dynamics can be written in a form
$$
  \dot x = -x(1-cu)+u
$$
and the Bellman function -- a solution of an optimal control problem $\tilde U(p,n,t)=\int_{T-t}^T (1-u(s))p(s)\mathrm ds$ with the starting point at $(p(t),n(t))=(p,n)$, can be present as $\tilde U(p,n,t)=nU(x,t)$.

\begin{figure}[t]
$$\includegraphics{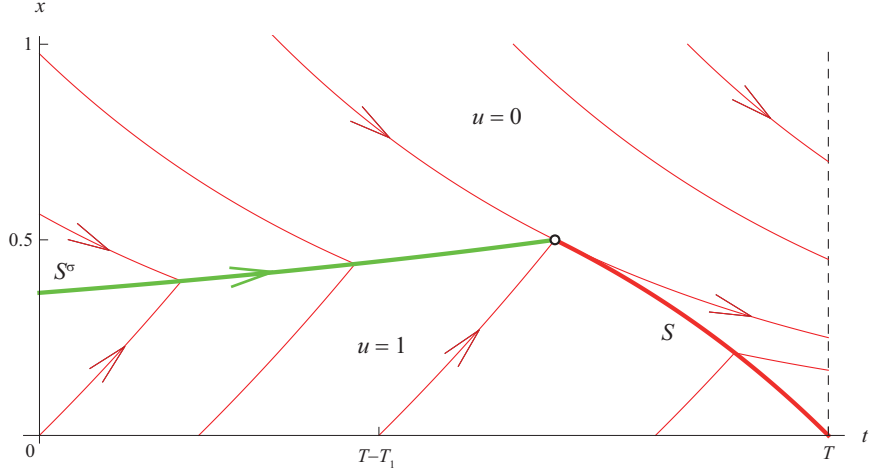}$$
\caption{Optimal collective behavior of the residents}
\label{f:2}
\end{figure}

The solution of the optimal control problem (\ref{e:dynamics0}-\ref{e:value0}) has been obtained before \cite{Akhmetzhanov2010} and the optimal behavioral pattern for $z=1.5$ and $T=2$ is shown on Fig.~\ref{f:2}. The region with $u=1$ is separated from the region with $u=0$ by a switching curve $S$ and a singular arc $S^\sigma$ such that
\begin{equation}
  \label{e:S}
S{:}\quad x=1-e^{-(T-t)}
\end{equation}
\begin{equation}
  \label{e:Ssigma}
 S^\sigma{:}\quad T-t = -\ln x+\frac2{xc}-\frac4c
\end{equation}
They are shown on Fig.~\ref{f:2} by thick red and green curves correspondingly. Along the singular arc $S^\sigma$ the consumer uses intermediate control $u=\hat u$:
\begin{equation}
\label{e:usigma}
\hat u=\frac{2x}{2+xc}
\end{equation}

One might identify a bang-bang control pattern for \emph{short seasons} $T\le T_1$ and a bang-singular-bang pattern for \emph{long seasons} $T>T_1$. The value $T_1$ is equal to
\begin{equation}
\label{e:T1}
T_1 = \frac{\ln(c+1)+(c-2)\ln2}{c-1}
\end{equation}
and it depends on the number of consumers present in the system.

The optimal value of the amount of offspring produced by individual can be computed using this solution. In the following we concentrate on the process of mutant invasion into population of consumers which uses the prescribed type of behavior given on Fig.~\ref{f:2}.

\subsection{Consumer-mutant-resource system}

Suppose that there is a subpopulation of consumers that acts as a mutant population. They maximize their own part of the fitness taking into account that the main population relates them as kin individuals.

Denote the fraction of the mutants with respect to the whole population of consumers by $\eps$ and variables describing a state of the mutant and resident populations by symbols with subindices ``$m$'' and ``$r$'' correspondingly. Then the number of mutants and residents will be $c_m=\eps c$ and $c_r=(1-\eps)c$ and the
dynamics of the system can be written in a form
\begin{equation}
  \label{e:dyn3}
  \dot p_r = -p_r+nu_r,\quad \dot p_m=-p_m+nu_m,\quad \dot n=-nc\left[(1-\eps)u_r+\eps u_m\right]
\end{equation}
similarly to (\ref{e:dynamics0}). The variable $u_m\in[0,1]$ defines a life-time decision pattern of the mutants. The control $u_r\in[0,1]$ is fixed and defined by the solution of the optimal control problem (\ref{e:dynamics0}-\ref{e:value0}).

The number of offspring for the next season is defined similarly to (\ref{e:value0}):
\begin{equation}
  \label{e:value3}
  J_r=\int_0^T \theta (1-u_r(t))c_rp_r(t)\>\mathrm dt,\ J_m=\int_0^T \theta (1-u_m(t))c_mp_m(t)\>\mathrm dt,\  J_n=\int_0^T \gamma n(t)\>\mathrm dt
\end{equation}
where the mutant chooses its control $u_m$ striving to maximize its criterion $J_m$.

We can see that the problem under consideration is described in terms of a two-step optimal control problem (or a hierarchical differential game): on the first step we define the optimal behavior of the residents, on the second step we identify the optimal response of the mutants to this strategy.

\section{Optimal free-riding}

Since $\theta$ and $\gamma$ are some constants, they can be omitted from the solution of the optimization problem $J_m\,\rightarrow\,\max\limits_{u_m}$. In this case the functional $J_m/(\theta c_m)$ can be taken instead of the functional $J_m$.

Let one introduce the Bellman function $\tilde U_m$ for the mutant population.
It provides a solution of the Hamilton-Jacobi-Bellman (HJB) equation
\begin{multline}
\label{e:HJB0}
\parder{\tilde U_m}{t}+\max\limits_{u_m}\left[ \parder{\tilde U_m}{p_r}(-p_r+nu_r)+\parder{\tilde U_m}{p_m}(-p_m+nu_m)-{}\right.\\\left.\parder{\tilde U_m}{n}nc((1-\eps)u_r+\eps u_m)+p_m(1-u_m)\right]=0
\end{multline}

Introducing new variables $x_r=p_r/n$ and $x_m=p_m/n$ and using a transformation of the Bellman function in the form $\tilde U_m(p_r,p_m,n,t)=nU_m(x_r,x_m,t)$, we can reduce the dimension of the problem by one using Noether's theorem. The modified HJB-equation (\ref{e:HJB0}) takes the following form
\begin{multline}
\label{e:HJB1}
\parder{U_m}{t} + \max\limits_{u_m}\Big\{ \parder{U_m}{x_r}\left[-x_r(1-c((1-\eps)u_r+\eps u_m))+u_r\right]+{}\\\parder{U_m}{x_m}\left[-x_m(1-c((1-\eps)u_r+\eps u_m))+u_m\right]-{}\\
U_mc((1-\eps)u_r+\eps u_m)+x_m(1-u_m)\Big\}=0
\end{multline}
Since the boundary conditions are defined at the terminal time it is convenient to construct the solution in backward time $\tau=T-t$. If we denote the components of the Bellman function as $\partial U_m/\partial x_r=\lambda_r$, $\partial U_m/\partial x_m=\lambda_m$ and $\partial U_m/\partial\tau = \nu$, equation (\ref{e:HJB1}) can be written as follows
\begin{multline}
\label{e:HJB2}
\cH\doteq-\nu + \max\limits_{u_m}\Big\{ \lambda_r\left[-x_r(1-c((1-\eps)u_r+\eps u_m))+u_r\right]+{}\\\lambda_m\left[-x_m(1-c((1-\eps)u_r+\eps u_m))+u_m\right]-{}\\
U_mc((1-\eps)u_r+\eps u_m)+x_m(1-u_m)\Big\}=0
\end{multline}
where the optimal control is defined as
\begin{equation*}
u_m = \Heav(\cA_m),\quad \cA_m=\partial\cH/\partial u_m=\lambda x_r\eps z+\lambda_m(1+x_m\eps c)-U_m\eps c-x_m\,.
\end{equation*}

One of the efficient ways to solve the HJB-equation is to use the method of characteristics (see \eg \cite{Melikyan1998}). The system of characteristics for equation (\ref{e:HJB2}) reads
\begin{equation}
\label{e:cs1}
  \begin{aligned}
    &x_r' = -\partial\cH/\partial\lambda_r = x_r(1-c((1-\eps)u_r+\eps u_m))-u_r,\\
    &x_m'= -\partial\cH/\partial\lambda_m = x_m(1-c((1-\eps)u_r+\eps u_m))-u_m,\\
    &\lambda_r'=\partial\cH/\partial x_r+\lambda_r\partial\cH/\partial U_m=-\lambda_r,\\
    &\lambda_m'=\partial\cH/\partial x_m+\lambda_m\partial\cH/\partial U_m=-\lambda_m+1-u_m,\\
    &\nu\,'=\nu\partial\cH/\partial U_m=-\nu c((1-\eps)u+\eps u_m),\\
    &U_m'=-\lambda_r\partial\cH/\partial\lambda_r-\lambda_m\partial\cH/\partial\lambda_m+\nu=\\
    &\hspace{2.5cm}-U_mc((1-\eps)u_r+\eps u_m)+x_m(1-u_m)
  \end{aligned}
\end{equation}
where the prime denotes differentiation with respect to backward time: $g'=dg/d\tau=-\dot g$. The terminal condition $U_m(x_r,x_m,T)=0$ gives that $\lambda_r(T)=\lambda_m(T)=0$. Then $\cA_m(T)<0$ and $u_m(T)=0$ as it could have been predicted before.

\subsection{First steps}

If we emit the characteristic field from the terminal surface $t=T$ with $u_r=u_m=0$ then
$$
x_r'=x_r,\quad x_m'=x_m,\quad \lambda_r'=-\lambda_r,\quad \lambda_m'=-\lambda_m+1,\quad U_m'=x_m\,,
$$
$$
\lambda_r(T)=\lambda_m(T)=0,\quad U_m(T)=0\,.
$$
We get the following equations for state and conjugate variables and for the Bellman function
$$
x_r=x_r(T)\mathrm e^\tau,\quad x_m=x_m(T)\mathrm e^\tau,\quad \lambda_r=0,\quad \lambda_m=1-\mathrm e^{-\tau},\quad U_m=x_m(1-\mathrm e^{-\tau})\,.
$$
\begin{figure}[t]
$$\includegraphics[scale=.9]{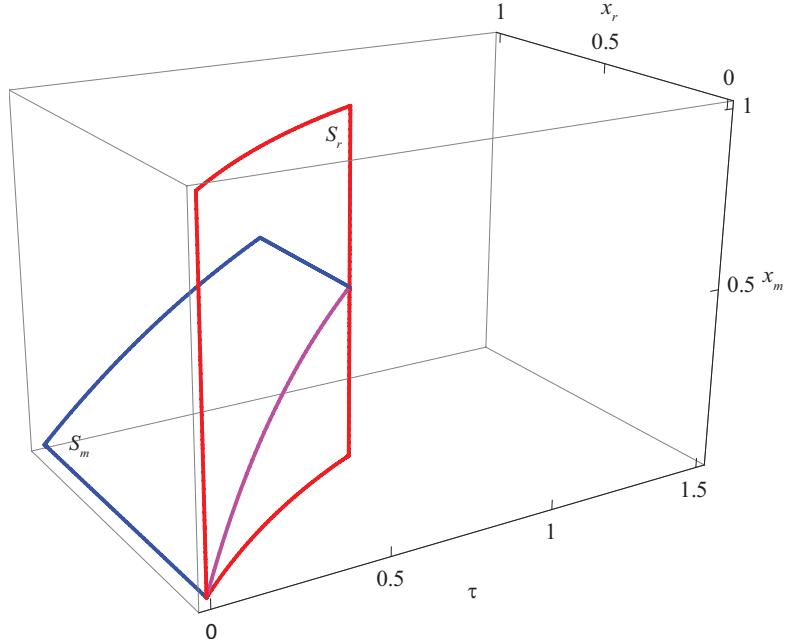}$$
\caption{First steps of the construction of the solution}
\label{f:3}
\end{figure}
From this solution we can see that there could exist a switching surface $S_m$:
\begin{equation}
\label{e:Sm}
S_m{:}\quad x_m=1-\mathrm e^{-(T-t)}
\end{equation}
such that $\cA_m=0$ on it and where the mutant is changing its control.
Equation (\ref{e:Sm}) is similar to (\ref{e:S}). But we should take into account the fact that there is also a hypersurface $S_r$, where the resident changes its control from $u_r=0$ to $u_r=1$ independently on the decision of the mutant. Hence it is important to define on which surface $S_r$ or $S_m$ the characteristic comes first, see Fig.~\ref{f:3}. Suppose that this is the surface $S_r$. Since the control $u_r$ changes its value on $S_r$, the HJB-equation (\ref{e:HJB2}) is also changing and, as a consequence, the conjugate variables $\nu$, $\lambda_r$ and $\lambda_m$ could possibly have a jump in their values. Let one denote the incoming characteristic field (in backward time) by ``$-$'' and the outcoming field by ``$+$''. Consider a point of intersection of the characteristic and the surface $S_r$ with the coordinates $(x_{r_1},x_{m_1},\tau_1)$. Then $x_{r_1}=1-\mathrm e^{-\tau_1}$ and the normal vector $\vartheta$ to the switching surface is written in the form
$$
\vartheta = \nabla S_r=(\partial S_r/\partial x_r,\partial S_r/\partial x_m,\partial S_r/\partial x_m)^T=(-1,0,1-x_{r_1})^T.
$$

From the incoming field we have the following information about the co-state
$$
\lambda_r^-=0,\quad \lambda_m^-=x_{r_1},\quad \nu^-=x_{m_1}\mathrm e^{-\tau_1}=x_{m_1}(1-x_{r_1})\,.
$$

Since the Bellman function is continuous on the surface $S_r$ which means
$$
U_m^+=U_m^-=U_m=x_{m_1}(1-\mathrm e^{-\tau_1})=x_{m_1}x_{r_1}\,.
$$
The gradient $\nabla U_m$ has a jump in the direction of the normal vector $\vartheta$: $\nabla U_m^+ = \nabla U_m^- + k\vartheta$. Here $k$ is an unknown scalar. Then
\begin{equation}
\label{e:lambda.plus}
\lambda^+ = -k,\quad \lambda_m^+ = x_{r_1},\quad \nu^+ = x_{m_1}(1-x_{r_1})+k(1-x_{r_1})
\end{equation}

If we suppose that the control of the mutant will be the same $u_m^+=0$ (in this case $\cA_m^+$ should be negative), the HJB-equation (\ref{e:HJB2}) has the form
\begin{equation}
\label{e:HJB2.plus}
-\nu^++\lambda^+[-x_{r_1}(1-(1-\eps)c)+1]-\lambda_m^+x_{m_1}(1-(1-\eps)c)-(1-\eps)zU_m+x_{m_1}=0
\end{equation}
By the substitution of the values from (\ref{e:lambda.plus}) to the equation (\ref{e:HJB2.plus}) we get
$$
k[-2(1-x_{r_1})-x_{r_1}(1-\eps)c]=0\,,
$$
which leads to the fact that $k=0$ and, actually, there is no jump in conjugate variables. They keep the same values as  (\ref{e:lambda.plus}) and $\cA_m^+=\cA_m^-$.

But let one suppose that the mutant reacts on the decision of the resident and also changes its control on $S_r$ from $u_m^-=0$ to $u_m^+=1$. This is fulfilled if the inequality $\cA_m^+>0$ holds.

The HJB-equation (\ref{e:HJB2}) has the form
$$
-\nu^++\lambda_r^+[-x_{r_1}(1-c)+1]+\lambda_m^+[-x_{m_1}(1-c)+1]-zU_m=0\,.
$$
Substitution of the values $\nu^+$, $\lambda_r^+$ and $\lambda^+_m$ from (\ref{e:lambda.plus}) gives
$$
k = \frac{x_{r_1}-x_{m_1}}{x_{r_1}z+(1-x_{r_1})}
$$
and
$$
\cA_m^+ = \lambda_r^+x_{r_1}\eps c+\lambda_m^+(x_{m_1}\eps c+1)-\eps cU_m-x_{m_1}=(x_{r_1}-x_{m_1})\frac{(1-\eps)x_{r_1}c+(1-x_{r_1})}{x_{r_1}c+(1-x_{r_1})}\,,
$$
which is positive when $x_{r_1}>x_{m_1}$. On Fig.~\ref{f:3} this corresponds to the points of the surface $S_r$ which are below the magenta line: $x_r=x_m=1-\mathrm e^{-\tau}$. For the optimal trajectories which go through such points $u_r(\tau_1+0)=u_m(\tau_1+0)=1$. One can show that there will be no more switches of the control. But if we consider a trajectory going from a point above the magenta line then $u_r(\tau_1+0)=1$ and $u_m(\tau_1+0)=0$ and there will be a switch of the control $u_m$ from zero to one (in backward time). After that there will be no more switches.\medskip

\begin{figure}[t]
$$\includegraphics[scale=.9]{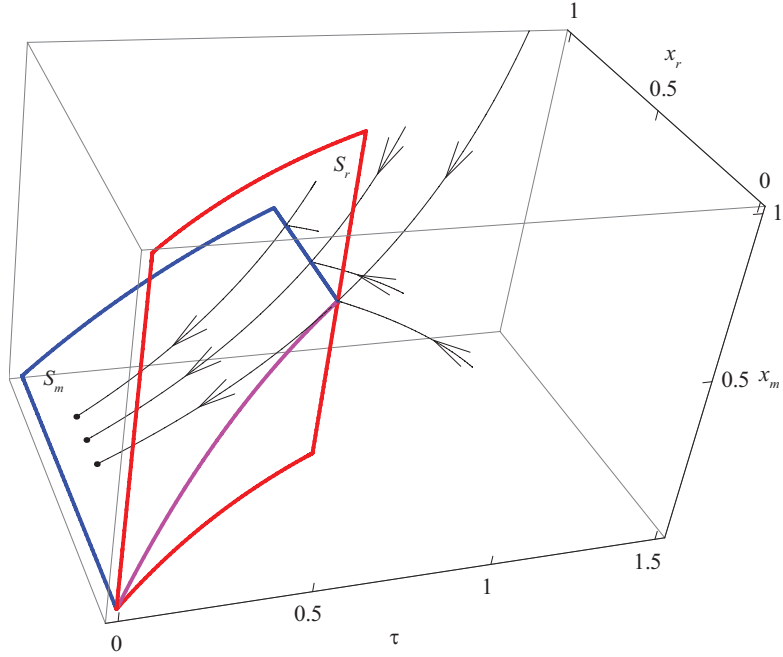}$$
\caption{Some family of the optimal trajectories emanated from the terminal surface}
\label{f:4}
\end{figure}

Now consider a trajectory emitted from the terminal surface which comes to the surface $S_m$ rather than to the surface $S_r$ at first. In this case the following situation as it is shown of Fig.~4 takes place. One might expect to have a singular arc $S_1^\sigma$ there.
Necessary conditions for its existence are the following
\begin{equation}
\label{e:sa1}
\cH = 0 = \cH_0 + \cA_mu_m,\quad \cH_0 = -\nu-\lambda x_r-\lambda_mx_m+x_m
\end{equation}
\begin{equation}
\label{e:sa2}
\cA_m=0=\lambda_rx_r\eps c + \lambda_m(x_m\eps c+1)-\eps cU_m-x_m
\end{equation}
\begin{equation}
\label{e:sa3}
\cA_m'=\{\cA_m\>\cH_0\}=0\doteq\cA_{m1}
\end{equation}
where the curled brackets denote the Poisson (Jacobi) brackets. If $\xi$ is a vector of state variables and $\psi$ is a vector of conjugate ones (in our case $\xi = (x_r,x_m,\tau)$ and $\psi = (\lambda_r,\lambda_m,\nu)$), then the Poisson brackets of two functions $F=F(\xi,\psi,U_m)$ and $G=G(\xi,\psi,U_m)$ are given by the formula
$$
\{F\,G\} = \langle F_\xi+\psi\,F_{U_m},\,G_\psi\rangle -\langle F_\psi,\,G_\xi+\psi\,G_{U_m}\rangle\,.
$$
Here $\langle \cdot,\cdot\rangle$ denotes the scalar product and \eg
$$
F_\psi=\partial F/\partial\psi=(\partial F/\partial\lambda_r,\partial F/\partial\lambda_m,\partial F/\partial\nu)^T.
$$

After some calculations the expression (\ref{e:sa3}) takes the form
\begin{equation}
\label{e:sa4}
\cA_{m1}=
\nu \eps c+x_m+\lambda_rx_r\eps c-(x_m+1)(1-\lambda_r)=0
\end{equation}

We can derive the variable $\nu$ from equation (\ref{e:sa1}) and substitute it to the last equation (\ref{e:sa4}). We get
$$
\cA_{m1} = x_m-1+\lambda_m=0\,.
$$
This leads to $\lambda_m=1-x_m$ and
$$
\lambda_r = \frac{x_m+\eps U_m+(1-x_m)(x_m\eps c+1)}{x_r\eps z}\,,
$$
which can be obtained from equation (\ref{e:sa2}).

To derive the singular control $u_m=\tilde u_m\in(0,1)$ along the singular arc one should write the second derivative
$$
\cA_m''=0=\{\{\cA_m\cH\}\cH\}=\{\{\cA_m\cH_0\}(\cH_0+\cA_m\tilde u_m)\}=\{\{\cA_m\cH_0\}\cH_0\}+\{\{\cA_m\cH_0\}\cA_m\}\tilde u_m\,.
$$
Then
\begin{equation}
\label{e:umsigma}
\tilde u_m = \frac{\{\{\cA_m\cH_0\}\cH_0\}}{\{\cA_m\{\cA_m\cH_0\}\}}=\frac{2x_m}{2+x_m\eps c}
\end{equation}
which has the same form as (\ref{e:usigma}).

The equation for the singular arc $S_1^\sigma$ can be obtained from dynamic equations (\ref{e:cs1}) by substitution $u_r=0$ and $u_m=\tilde u_m$ from (\ref{e:umsigma}):
$$
x_m'=-\frac{x_m^2\eps c}{2+x_m\eps c},\quad x_m(\tau=\ln2)=1/2\,.
$$
Finally, we have the analogous expression to (\ref{e:Ssigma})
\begin{equation}
  \label{e:Smsigma}
  S_1^\sigma{:}\quad T-t= -\ln x_m+\frac2{x_m\eps c}-\frac{4}{\eps c}
\end{equation}
for $\eps\ne0$. If $\eps=0$ the surface $S_m$ is a hyperplane $x_m=1/2$.

\begin{figure}[t]
$$\includegraphics[scale=1]{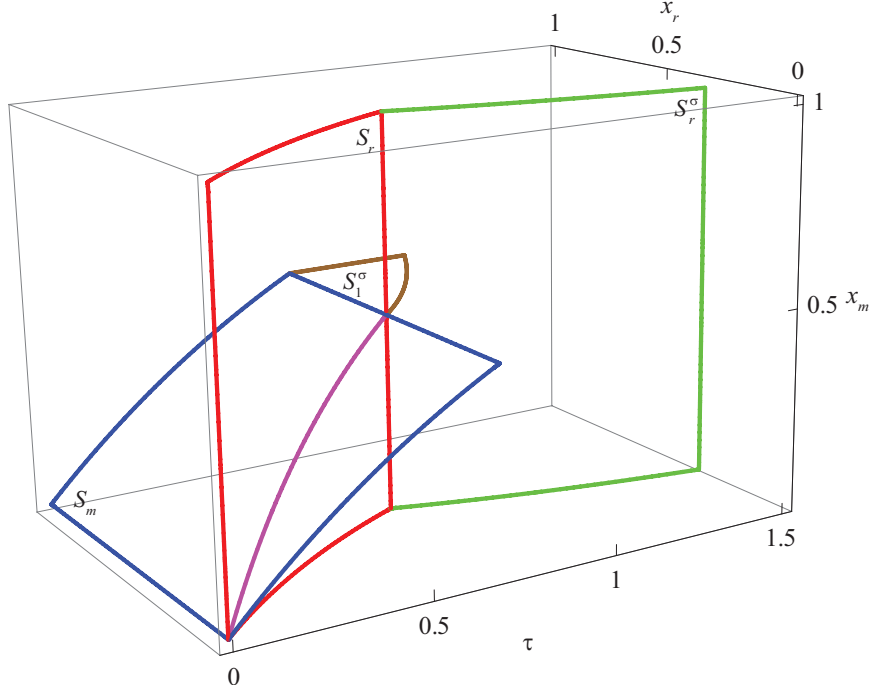}$$
\caption{Construction of the part of a singular arc $S_1^\sigma$}
\label{f:5}
\end{figure}

After these steps we have the structure of the solution shown on Fig.~\ref{f:5}.

\subsection{Optimal motion along the surface $S_r^\sigma$}

Let us consider the surface $S_r^\sigma$ which is shown on Fig.~\ref{f:5} by green color and separates the domain $u_r=1$ from $u_r=0$. This leads to a chattering regime with the motion along this surface. The solution of the dynamic equations can be understood in Fillipov sense.

Suppose that the hypersurface $S_r^\sigma$ is also divided into two regions: the region where the mutant uses $u_m=0$ and the one where $u_m=1$. They are separated by a switching curve $\hat S$ and by a singular arc $\hat S^\sigma$ which completely belong to the surface $S_r^\sigma$. Along this surface the resident uses an intermediate control $u_r=\hat u_r$ resulting from the chattering regime with simultaneous switches from one bang control $u_r=1$ to another $u_r=0$ and vice versa. We suppose that the trajectory can be forced to stay on the surface $S_r^\sigma$  by the resident independently on the action of the mutant. This means that if we derive the control $\hat u_r$ from the dynamic equation
$$
x_r' = -\frac{x_r^2c}{2+x_rc}=x_r(1-c((1-\eps)\hat u_r+\eps u_m))-\hat u_r
$$
as
\begin{equation}
\label{e:hatu}
\hat u_r = \frac{2x_r(1+x_rc)}{(1+(1-\eps)x_rc)(2+x_rc)}-\frac{x_r\eps cu_m}{1+(1-\eps)x_rc}
\end{equation}
Then $\hat u_r\in (0,1)$ for all points belonging $S_r^\sigma$ and for all possible values $u_m\in[0,1]$.

To identify for which parameters of the model this is possible, we may notice that $\hat u_r$ as a function of $u_m$ is linear and decreasing. Moreover
$$
u_r(u_m=0) = \frac{2x_r(1+x_rc)}{(1+(1-\eps)x_rc)(2+x_rc)}\le 2x_r\frac{1+x_rc}{2+x_rc}\le1\,.
$$
since $x_r\le1/2$. Therefore only the condition $u_r\ge0$ could be violated for some values of $u_m$. To define the limiting value $\tilde u_m$ for which $u_r(\tilde u_m<u_m\le1)<0$ one can use the following condition: $u_r(u_m=\tilde u_m)=0$. This gives
$$
\tilde u_m= \frac2{\eps c}\frac{1+x_rc}{2+x_rc}\,.
$$
If such value $\tilde u_m$ is outside of the interval $(0,1)$ then the condition $u_r\in(0,1)$ holds for any $x_r$ belonging to $S_r^\sigma$. This occurs if
\begin{equation}
\label{e:epsbound}
\eps<1/c
\end{equation}

In this paper we consider only the values of $\eps$ satisfying (\ref{e:epsbound}). This has a biological explanation since for sufficiently large $\eps$ the resident should react to the behavior of the mutants who does significantly significantly affect the dynamics of the system. A fixed \emph{a priory} strategy of the resident does not make sense in that case.\medskip

If we consider a field belonging to the surface $S_r^\sigma$, the gradient of the restriction of $U_m$ to that manifold is defined only in the co-tangent bundle. A safe representative requires a term $k n_s$ be added to the adjoint equations of the characteristic system, where $n_s$ is the normal to $S_r^\sigma$. The constant $k$ should be chosen to keep the adjoint tangent to it. But we can notice that the surface $S_r^\sigma$ does not depend on $x_m$-coordinate. Since only the dynamics $\lambda_m$ plays an important role for us and the corresponding term is equal zero, this notion can be neglected.

The control $\hat u_r = \hat u_r(x_r,x_m,\tau,u_m)$ is defined in feedback form, \eg depends on time and a state of the system. The corresponding Hamiltonian (\ref{e:HJB2}) is changed to
\begin{equation}
\label{e:hatHamiltonian}
\hat\cH = \cH(x_r,x_m,U_m,\lambda,\lambda_m,\nu,\hat u_r(x_r,x_m,\tau,u_m),u_m)\,.
\end{equation}
The coefficient multiplying the control $u_m$ is also changed to
\begin{equation}
  \label{e:hatAm}
  \hat \cA_m = \parder{\hat\cH}{u_m}= \frac{\lambda_m(1+x_r(1-\eps)c+x_m\eps c)-\eps cU_m}{1+(1-\eps)x_rc}-x_m\,.
\end{equation}

In this case the switching surface $\hat S$ can be defined by the condition $\hat \cA_m=0$. The singular arc $\hat S^\sigma$ -- by the following conditions
\begin{equation*}
\label{e:sa_nc_hat}
\hat\cH = 0,\quad \hat \cA_m=0,\quad \hat\cA_m'=\{\hat\cA_m\hat\cH\}=0\,.
\end{equation*}
The intermediate control $\hat u_m$ can be obtained from the second derivative
$$
\hat\cA_m''=\{\{\hat\cA_m\hat\cH\}\hat\cH\}=0\,.
$$
We can write analytical expressions for $\hat S$ and $\hat S^\sigma$ but they look quite complicated. To make things simpler, let us consider first a particular case of vanishingly small values of $\eps$ and study the optimal behavioral pattern.

\subsection{Particular case of a vanishingly small population of mutants}

We have $\eps\cong0$ and the chattering regime of the resident along the surface $S_r^\sigma$ results in $u_r=\hat u_r$ coinciding with (\ref{e:usigma}):
$$
\hat u_r = \frac{2x_r}{2+x_rc}
$$
which does not depend on the action of the mutant.
In addition equations (\ref{e:hatHamiltonian}) and (\ref{e:hatAm}) take the following form
\begin{equation}
\label{e:hatHamiltonian2}
\hat\cH = -\nu + \frac{\lambda_r x_r^2c}{2+x_rc}+\lambda_m\left(-x_m\frac{2-x_rc}{2+x_rc}+u_m\right)-U_m\frac{2x_rc}{2+x_rc}+x_m(1-u_m)\,,
\end{equation}
\begin{equation}
\label{e:hatAm2}
\hat\cA_m=\lambda_m-x_m\,.
\end{equation}

If the trajectory goes from the point $x_m^\sigma\doteq x_m(\ln2)>1/2$ then $u_m=0$ and the system of characteristics for the Hamiltonian (\ref{e:hatHamiltonian}) is written in the form
\begin{equation*}
\label{e:hatSC}
x_r'=-\frac{x_r^2c}{2+x_rc},\quad x_m'=x_m\frac{2-x_rc}{2+x_rc},\quad \lambda_m=-\lambda_m+1\,,
\end{equation*}
$$
\nu\,'=-\nu\frac{2x_rc}{2+x_rc},\quad U_m'=-U_m\frac{2x_rc}{2+x_rc}+x_m\,,
$$
with boundary conditions
\begin{equation*}
\label{e:hatBC}
\tau=\ln2,\quad x_r=1/2,\quad x_m=x_m^\sigma,\quad \lambda_m=1/2,\quad \nu=x_m^\sigma/2,\quad U_m=x_m^2/2\,.
\end{equation*}
Then $\lambda_m=1-\mathrm e^{-\tau}$ and the switching curve $\hat S$ has the form
$$
x_m=1-\mathrm e^{-\tau},\quad
$$
in addition to $\tau = -\ln x_r + 2/(x_rc)-4/c$. Thus $\hat S=S_m\cap S_r^\sigma$.

The switching curve $\hat S$ ends at the point with coordinates $(x_{r_2},x_{m_2},\tau_2)$ where the characteristics become tangent to it and the singular arc $\hat S^\sigma$ appears. Before the determination of the coordinates of this point let one define the singular arc $\hat S^\sigma$. From equations (\ref{e:hatHamiltonian2}) and (\ref{e:hatAm2}) we get
\begin{equation*}
\label{e:sa11}
\nu = \frac{\lambda_rx_r^2c}{2+x_rc}-\lambda_mx_m\frac{2-x_rc}{2+x_rc}-U_m\frac{2x_rc}{2+x_rc}+x_m,\quad \lambda_m=x_m
\end{equation*}
along the singular arc.
Substitution of (\ref{e:sa11}) into equation $\hat\cA_m'=0$ gives
\begin{equation*}
\label{e:hatxm}
  x_m=\frac{2+x_rc}4\,.
\end{equation*}
Alongside, the intermediate control $\hat u_m$ can be derived from $\hat\cA_m''=0$ and it is equal to
$$
\hat u_m=\frac1{2+x_rc}\,,
$$
which is positive and belongs to the segment between zero and one.

\begin{figure}[t]
$$\includegraphics[scale=.9]{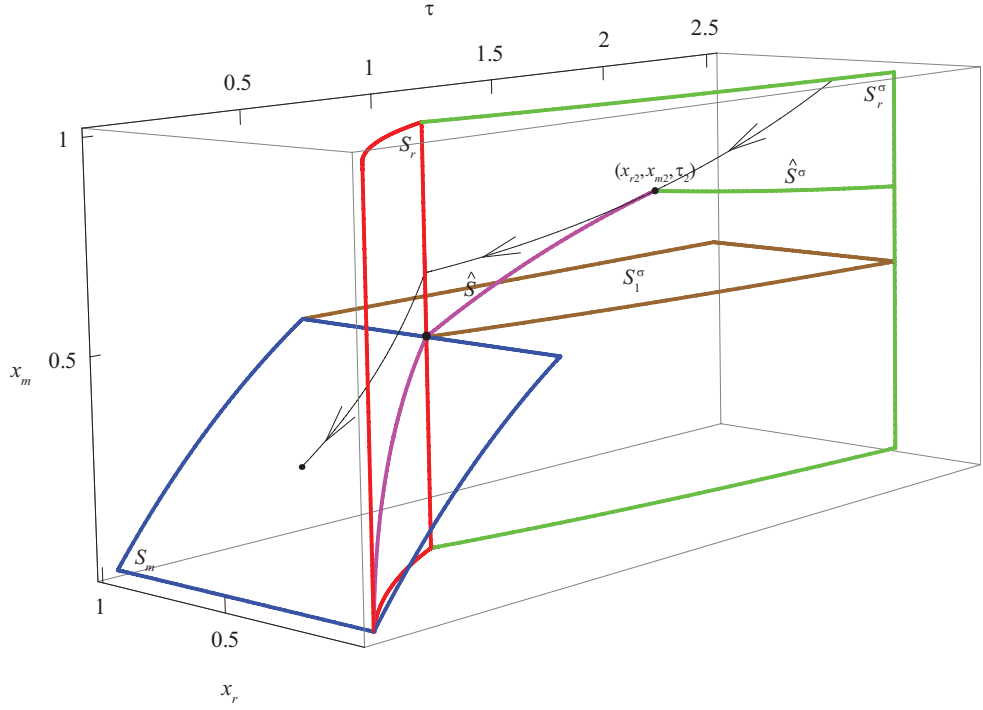}$$
\caption{Optimal behavior inside the surface $S_r^\sigma$}
\label{f:6}
\end{figure}

We see that the coordinates $x_{r_2}$, $x_{m_2}$ and $\tau_2$ can be defined through the following equations
\begin{equation*}
\label{e:pointB}
x_{m_2}=\frac{2+x_{r_2}c}4=1-\mathrm e^{-\tau_2},\quad \tau_2=-\ln x_{r_2}+\frac2{x_{r_2}c}-\frac4c\,,
\end{equation*}
which comes from the fact that the point $(x_{r_2},x_{m_2},\tau_2)$ belongs to $\hat S^\sigma$ and it is located on the intersection of the curves $\hat S^\sigma$ and $\hat S$. The result is illustrated on Fig.~\ref{f:6}.

We can show also that the surface $S_1^\sigma$ can be extended further with comparison to the situation on Fig.~\ref{f:5}. Indeed, the following conditions are fulfilled for the region with $u_r=0$:
$$
\cH\Big|_{u_r=0}=-\nu-\lambda_rx_r-\lambda_mx_m+x_m=0,\quad \cA_m=\lambda_m-x_m=0,\quad \cA_m'=0\,.
$$
Therefore
$$
\nu=-\lambda_rx_r-\lambda_mx_m+x_m,\quad \lambda_m=x_m\,,
$$
and from the condition $\cA_m'=0$: $-1+2x_m=0$, which gives $x_m=1/2$.

\begin{figure}[t]
$$\includegraphics[scale=.9]{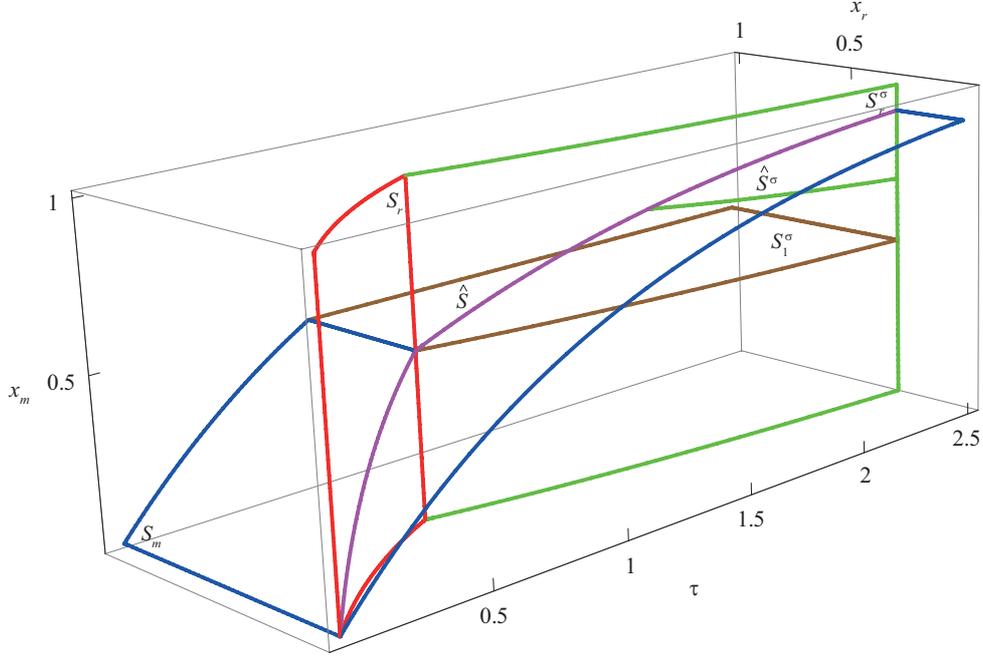}$$
\caption{Optimal behavioral pattern for $c=3$}
\label{f:7}
\end{figure}

Consider now the region with $x_r$ smaller that the ones on the green surface $S_r^\sigma$ (see Fig.~\ref{f:6}). There is a switching surface which extends the surface $S_m$ and it is defined by the same equation (\ref{e:Sm}). But there could exist a singular arc $S_2^\sigma$ starting from some points of $S_m$. To check this we have to write the following conditions
\begin{equation}
\label{e:nc1}
\cH\Big|_{u_r=1} = -\nu-\lambda_r(x_r(1-c)-1)-\lambda_mx_m(1-c)-zU_m+x_m=0\,,
\end{equation}
\begin{equation}
\label{e:nc2}
\cA_m=\lambda_m-x_m=0,\quad \cA_m'=0\,,
\end{equation}
which give a possible candidate for a singular arc
$$
S_2^\sigma{:}\quad x_m=\frac1{2-c}\,.
$$
We see that its appearance is possible only for $c<2$. In addition, the motion along this surface occurs with control $u_m=(1-c)/(2-c)$ which also gives the restriction on the parameter $c$ that $c\le 1$. For $c>1$ the structure of the solution in the domain below the surface $S_r^\sigma$ is simpler and consists only of the switching surface $S_m$, see Fig.~\ref{f:7}.

\subsection{Computation of the value functions in case of $\eps=0$}
\label{s:valuefunctions}

\begin{figure}[p]
$$\hspace{-.5cm}\includegraphics[scale=.95]{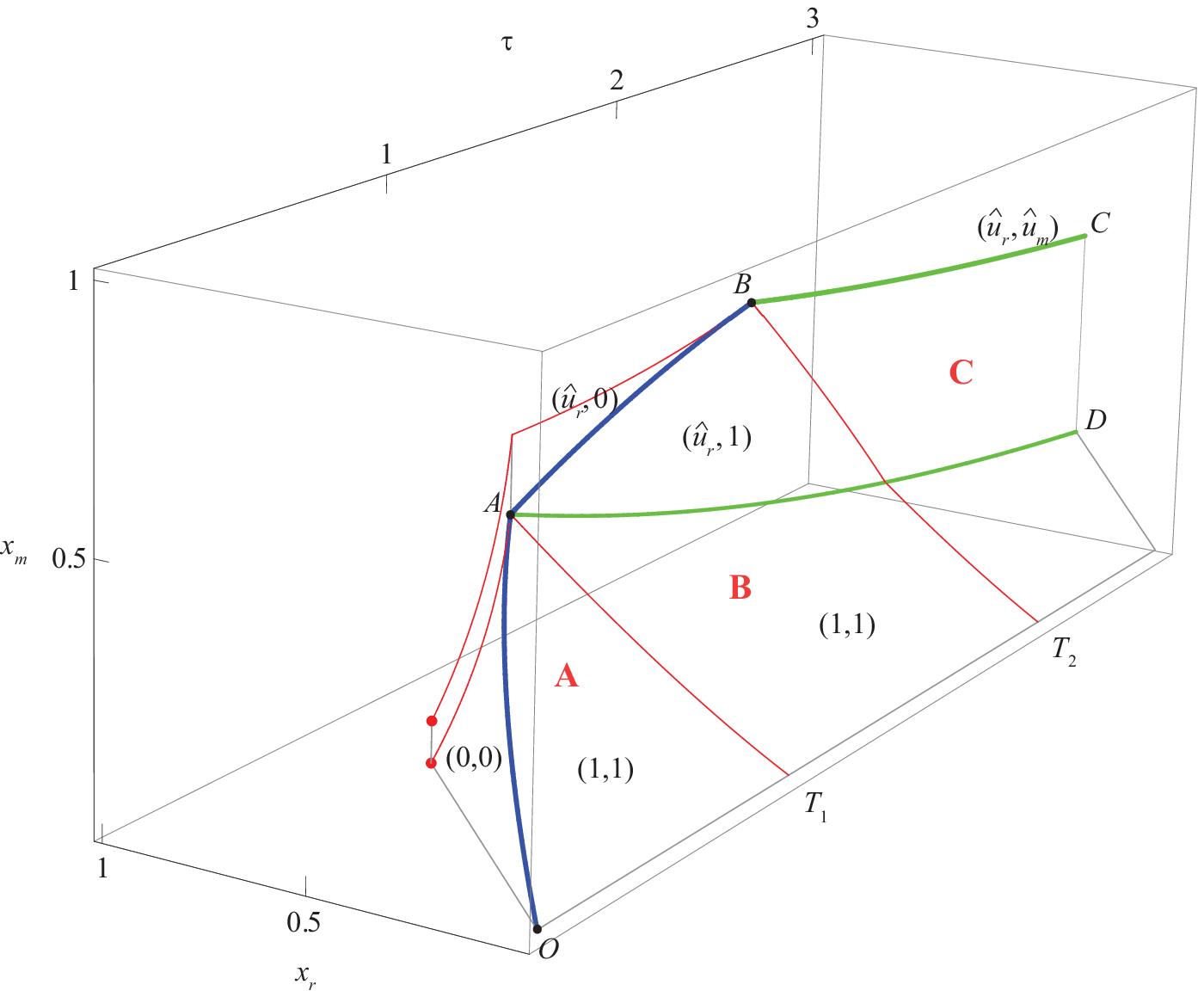}$$
\caption{The reduced optimal pattern for trajectories coming from zero initial conditions $x_r(0)=x_m(0)=0$ and $c=3$}
\label{f:8}
$$\includegraphics[scale=.9]{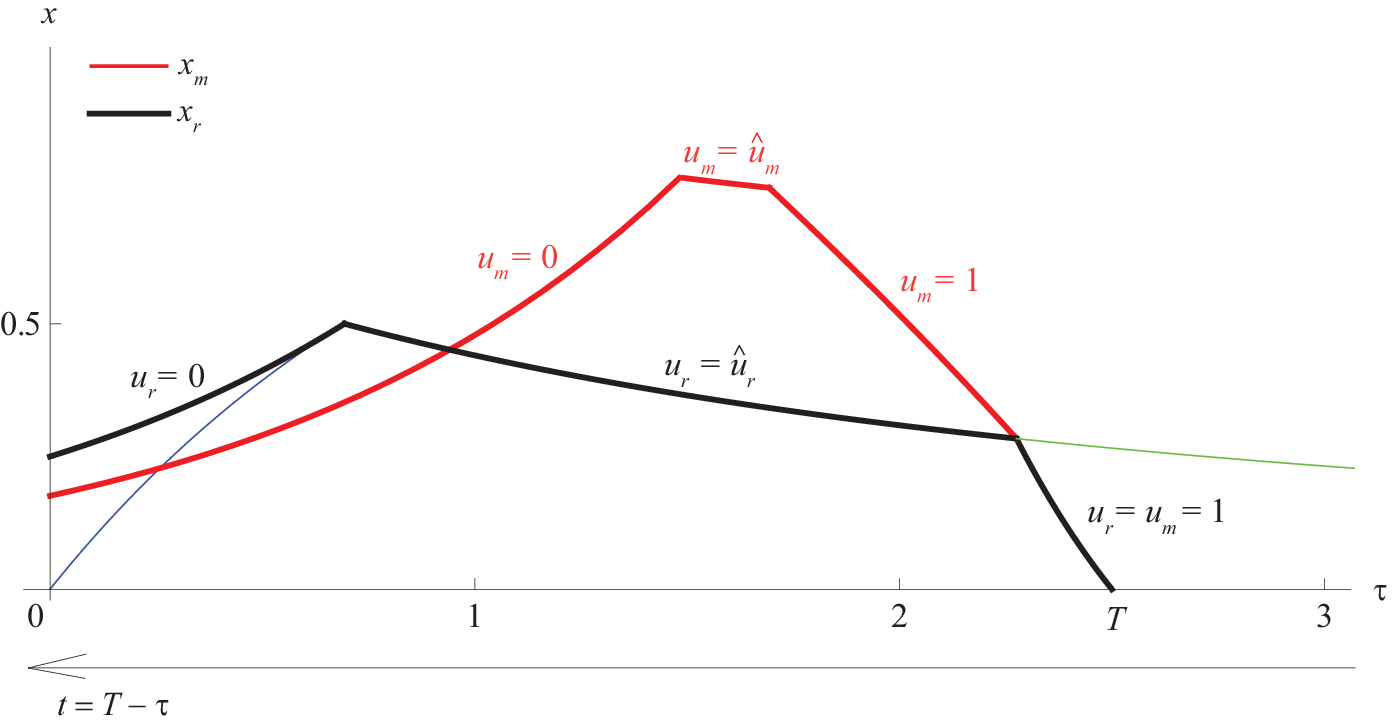}$$
\caption{Optimal free-riding of the mutant}
\label{f:9}
\end{figure}

Without loss of generality we can assume that at the beginning of each season the average energy of the population of consumers is zero: $x_r(0)=x_m(0)=0$. Therefore we should take into account only the trajectories coming from these zero initial conditions. The phase space is reduced in this case to the one shown on Fig.~\ref{f:8}. One can see that there are three different regions depending on the length of the season $T$. If it is short enough $T\le T_1$ (where the value $T_1$ has been defined in (\ref{e:T1})), then the behavior of the mutant coincides with the behavior of the resident and the main population can not be invaded: the amount of offspring produced by the mutant is the same as produced by the resident. If the length of the season is larger than $T>T_1$ there is a period of the life-time of the resident when it applies the intermediate strategy and spares some amount of the resource for its future use. The mutant is able to use this fact and there exists a strategy of the mutant that guarantees better result for it.

Let us introduce the analogue of the value function $\tilde U_m$ for the resident and denote it as $\tilde U_r$:
$$
\tilde U_r(p_r,p_m,n,t) = \int_{T-t}^T p_r(s)(1-u_r(s))\>\mathrm dt\,.
$$
The value $\tilde U(0,0,n(0),T)$ represents the amount of eggs laid by the resident during the season of length $T$. Its value depends on the state of the system and the following transformation can be done
$$
\tilde U_r(p_r,p_m,n,t) = n U_r(x_r,x_m,t)\,.
$$
In the following we omit some parameters and write the value function in the simplified form $U_r(T)\doteq U_r(0,0,T)$ where the initial conditions $x_r(0)=x_m(0)=0$ have been taken into account.

In the region $\mathbf A$  (see Fig.~\ref{f:8}) the value functions for both populations of mutants and residents are equal to each other
$$
U_m(T)=U_r(T)=x_1\mathrm e^{-c(T-\tau_1)}\,.
$$
Here the value $\tau_1$ can be defined from the intersection of the trajectory and the switching curve $S_r\cap S_m$:
$$
1-\mathrm e^{-\tau_1}=\frac{\mathrm e^{(c-1)(T-\tau_1)}-1}{c-1}\,.
$$

To obtain the value functions in the regions $\mathbf B$ and $\mathbf C$ let one solve the system of characteristics (\ref{e:hatSC}) in case when the characteristics move along the surface $S_r^\sigma$ and $u_m=1$. This leads to the following characteristic equations for the Hamiltonian (\ref{e:hatHamiltonian2}):
$$
x_r'=-\frac{x_r^2c}{2+x_rc},\quad x_m'=x_m\frac{2-x_rc}{2+x_rc}-1,\quad U_m'=-U_m\frac{x_r^2c}{2+x_rc}\,.
$$
We can rewrite them in the form
$$
\frac{dx_m}{dx_r}=\frac{2+x_rc-x_m(2-x_rc)}{x_r^2c},\quad \frac{dU_m}{dx_r}=\frac{2U_m}{x_r}\,,
$$
and consequently
\begin{equation}
\label{e:dynSsigma}
x_m=C_1x_r^2\mathrm e^{\tau}+x_rz+1,\quad U_m=C_2x_r^2,\quad C_1,C_2=\mathrm{const}\,.
\end{equation}
where $C_1$ and $C_2$ are defined from the boundary conditions while the equation (\ref{e:Ssigma}) is also fulfilled.

Along the singular arc $\hat S^\sigma$ the mutant is using the intermediate strategy (\ref{e:umsigma}). In this case
$$
U_m'=-U_mzu_r+x_m(1-u_m)=-U_m\frac{2x_rc}{2+x_rc}+\frac{1+x_rc}4\,.
$$
Since $x_r'=-\frac{x_r^2c}{2+x_rc}$ we have
$$
\frac{dU_m}{dx_r}=\frac{2U_m}{x_r}-\frac{(1+x_rc)(2+x_rc)}{4x_r^2c}\,.
$$
Then
\begin{equation}
\label{e:UmC3}
U_m=C_3x_r^2+\frac{4+3x_rc(3+2x_rc)}{24x_rc},\quad C_3=\mathrm{const}\,.
\end{equation}

We undertake now to compute the limiting season length $T_2$ that separates the region $\mathbf B$ from the region $\mathbf C$. The coordinates of the point $B$ have been obtained before and satisfy the equations (\ref{e:pointB}). To define the coordinates of the point $(x_{r_2}^\sigma,x_{m_2}^\sigma,\tau_2^\sigma)$ of intersection of the optimal trajectory with the curve $AD$ let use the dynamics of the motion along the surface $S_r^\sigma$ with $u_r=\hat u_r$ and $u_m=1$ (\ref{e:dynSsigma}):
$$
x_m=C_1x_r^2\mathrm e^\tau+x_rz+1,\quad C_1=\mathrm{const}\,,
$$
where $C_1$ should be chosen such that
$$
x_{m_2}=C_1x_{r_2}^2\mathrm e^{\tau_2}+x_{r_2}c+1,\quad x_{m_2}=\frac{2+x_{r_2}c}4=1-\mathrm e^{-\tau_2}\,.
$$
Then
\begin{equation*}
\label{e:C1}
C_1=\frac{(x_{r_2}c-2)(3x_{r_2}c+2)}{16x_{r_2}^2}\,.
\end{equation*}
After that the coordinates $x_{r_2}^\sigma$, $x_{m_2}^\sigma$ and $\tau_2^\sigma$ can be defined from the following conditions
$$
x_{m_2}^\sigma=x_2^\sigma=C_1(x_{r_2}^\sigma)^2\mathrm e^{\tau_2^\sigma}+x_{r_2}^\sigma c+1,\quad \tau_2^\sigma=-\ln x_{r_2}^\sigma+\frac2{x_{r_2}^\sigma c}-\frac4c\,.
$$
The boundary value $T_2$ can be obtained as
$$
T_2 = \tau_2^\sigma +\frac{\ln(x_{r_2}^\sigma(c-1)+1)}{c-1}\,.
$$

Now compute the value functions $U_r(T)$ and $U_m(T)$  for the region $\mathbf B$ ($T_1<T\le T_2$), where only the mutant keeps the bang-bang type of the control. For the resident population we have
\begin{equation}
\label{e:UT}
U_r(T) = U_{r_2}\mathrm e^{-c(T-\tau_2)},\quad U_{r_2} = x_{r_2}(1-x_{r_2})+\frac{1-2x_{r_2}}c\,.
\end{equation}
where the point with coordinates $(x_{r_2},x_{r_2},\tau_2)$ defines the intersection of the trajectory and surface $S_r^\sigma$:
\begin{equation}
\label{e:x2tau2}
\tau_2 = -\ln x_{r_2}+\frac2{x_{r_2}c}-\frac4c,\quad x_{r_2}=\frac{\mathrm e^{(c-1)(T-\tau_2)}-1}{c-1}\,.
\end{equation}
For the mutant population the value function $U_m$ in the region with $u=\hat u$ and $u_m=1$ satisfies the equation coming from (\ref{e:dynSsigma}):
\begin{equation}
\label{e:Um1}
U_m^{(\hat u,1)}=x_{m_1}^2\left(\!\frac{x_r}{x_{r_1}}\!\right)^2\,,
\end{equation}
where $(x_{r_1},x_{m_1},\tau_1)$ is a point of the intersection of the trajectory with the curve $AB$ (see Fig.~\ref{f:8}). Using (\ref{e:Um1}) and notations of (\ref{e:x2tau2}), we can write
$$
U_m(T) = U_{m_2}\mathrm e^{-c(T-\tau_2)},\quad U_{m_2}=x_{m_1}^2\left(\!\frac{x_{r_2}}{x_{r_1}}\!\right)^2\,,
$$
which is analogous to (\ref{e:UT}).

\begin{figure}[t]
$$\hspace{-.75cm}\includegraphics[scale=.75]{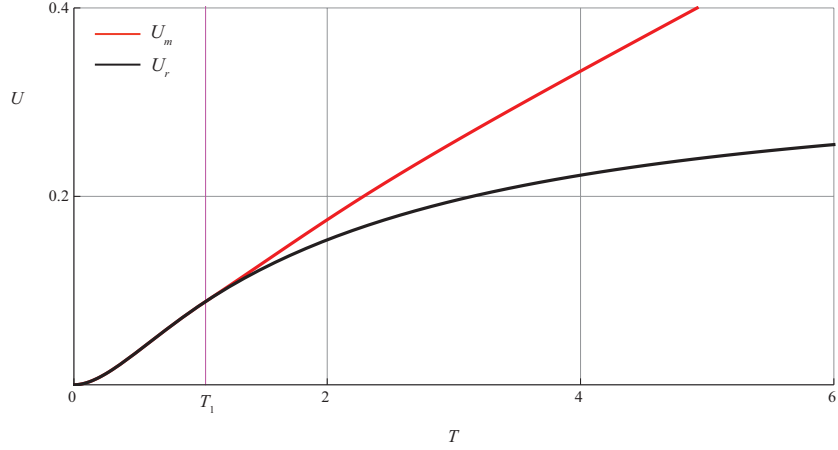}$$
\caption{Difference in value functions for the resident and the mutant ($c=3$)}
\label{f:10}
\end{figure}

For the region $\mathbf C$ the value function for the resident has the same form as in (\ref{e:UT}) but it has a different form for the mutant. Suppose that the optimal trajectory is coming to the surface $S^\sigma$ at point with coordinates $(\tilde x_{r_2},\tilde x_{m_2},\tilde \tau_2)$. Then the Bellman function at this point is equal
$$
\tilde U_{m_2} = \tilde x_{r_2}^2\left(\!\frac{c^2}{16}-\frac{4+3\tilde x_{r_2}c}{24\tilde x_{r_2}^3c}\right)+\frac{3\tilde x_{r_2}c(2\tilde x_{r_2}c+3)+4}{24\tilde x_{r_2}c}\,,
$$
which is written using (\ref{e:UmC3}) with definition of the constant $C_3$ from the given boundary conditions.

When the optimal trajectory moving along the surface $S^\sigma$ intersects the curve $AD$ at some point with coordinates $(\tilde x_{r_2}^\sigma,\tilde x_{m_2}^\sigma,\tilde \tau_2^\sigma)$ (see Fig.~\ref{f:8}) the Bellman function can be expressed as follows
$$
\tilde U_{m_2}^\sigma = \tilde U_{m_2}\left(\!\frac{\tilde x_{r_2}^\sigma}{\tilde x_{r_2}}\!\right)\,.
$$
Then
$$
U_m(T) = \tilde U_{m_2}^\sigma \mathrm e^{-c(T-\tau_2^\sigma)}\,.
$$

The difference in the values functions (amount of offspring per mature individual) of the mutant and optimally behaving resident is shown on Fig.~\ref{f:10}. In way we can derive the expressions for the number of offspring produced by the resource population during the season.

\subsection{Generalization on the small enough but non-zero values of $\eps$}

In this section we consider a case of non-zero $\eps$ but such that the condition (\ref{e:epsbound}) remains fulfilled. This means that the trajectory coming to the singular surface $S_r^\sigma$ does not cross it but moves along it due to the chattering regime applied by the resident (\ref{e:hatu}).

In this case the phase space can be also divided in two regions: the points with $x_r$ smaller or larger than the ones on $S_r^\sigma$. In each of this region the structure of the solution has similar properties as in the case considered above when $\eps$ is arbitrary small. Inside the surface $S_r^\sigma$ the optimal behavior has also a similar to a previous case structure.

In the region with values $x_r$ larger than the ones on the surface $S_r^\sigma$ there is a part of the switching surface $S_m$ and a singular arc $S_1^\sigma$ where the mutant uses an intermediate strategy. The surface $S_1^\sigma$ can be defined through the expression (\ref{e:Smsigma}). In the other region we also have a part of $S_m$ and a singular arc $S_2^\sigma$ which is different from $S_1^\sigma$ and could not exist for some values of the parameters of the problem $c$ and $\eps$.

To identify the values for which the surface $S_2^\sigma$ is a part of the solution let us write necessary conditions similarly to (\ref{e:nc1}-\ref{e:nc2}):
$$
\cH\Big|_{u_r=1}=0,\quad \cA_m=0,\quad \cA_m'=\{\cA_m\cH\}=0\,.
$$
Using these equations we are able to obtain the values of $\lambda_r$, $\lambda_m$ and $\nu$ on the surface $S_2^\sigma$ and substitute them into the second derivative $\cA_m''=\{\{\cA_m\cH\}\cH\}=0$ to derive the expression for the singular control applying by the mutant on this surface:
\begin{equation}
\label{e:um3}
u_m=\frac{2x_m-(1-\eps)c(1+x_m)}{2-(1-\eps)c+x_m\eps c}\,.
\end{equation}
There are several conditions which should be necessarily satisfied. First of all, the control (\ref{e:um3}) should be between zero and one
\begin{equation}
\label{e:um01}
0\le\frac{2x_m-(1-\eps)c(1+x_m)}{2-(1-\eps)c+x_m\eps c}\le1
\end{equation}
Second of all, the Kelley condition should be also fulfilled \cite[p.~200]{Melikyan1998}:
$$
\parder{}{u_m}\frac{d^2}{dt^2}\parder{\cH}{u_m}=\{\cA_m\{\cA_m\cH\}\}\le0\,.
$$
This leads to the inequality
\begin{equation}
\label{e:kel01}
2-(1-\eps)+x_m\eps c\ge0
\end{equation}

In particular, both conditions (\ref{e:um01}) and (\ref{e:kel01}) give that $x_m\le2/(2-c)$.

\begin{figure}[t]
$$\hspace{1.cm}\includegraphics[scale=.9]{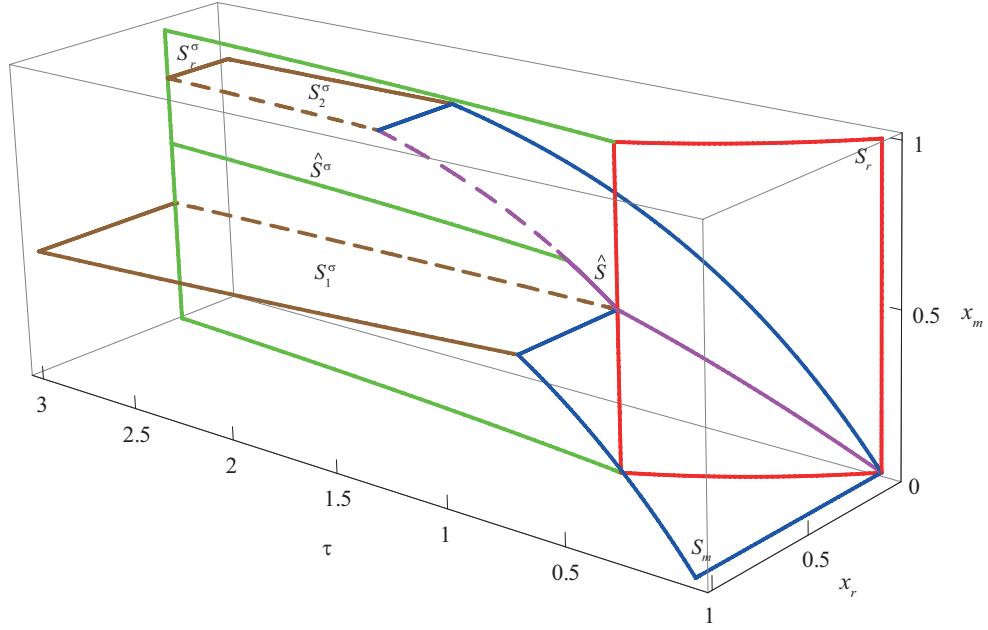}$$
\caption{Structure of the optimal behavioral pattern for $c=1.25$ and $\eps=0.35$}
\label{f:11}
\end{figure}

To construct the singular arc $S_2^\sigma$ we should substitute the singular control $u_m$ from (\ref{e:um01}) and $u_r=1$ into dynamics (\ref{e:cs1}):
$$
x_m' = x_m(1-c((1-\eps)u_r+\eps u_m))-u_m
$$
with boundary conditions obtained from the tangency condition of the optimal trajectory coming from the domain $u_m=u_r=1$ on the switching surface $S_m$:
$$
x_{m}\Bigl(-\ln\bigl(1-\frac1{2-c(1-\eps)}\bigr)\Bigr)=\frac1{2-c(1-\eps)}
$$
Such tangency takes place only if
$$
0\le\frac1{2-c(1-\eps)}\le1
$$
that comes from the condition that a singular surface $S_m$ exists only for $0\le x_m<1$. This gives the following inequality
$$
1-c(1-\eps)\ge0
$$
for the existence of the surface $S_2^\sigma$. One can check that the inequalities (\ref{e:um01}-\ref{e:kel01}) are fulfilled as well. The result of the construction of the optimal pattern for some particular case is shown on Fig.~\ref{f:11}.

\section{Long-term evolution of the system}

\cite{Akhmetzhanov2010} introduced model (\ref{e:dynamics0}) as the intra-seasonal part of a more complex multi-seasonal population dynamics model in which consumer and resources survive during one season only. They considered that the (immature) offspring produced by the consumers and ressources through some season $i$ and defined by equations (\ref{e:value0}) would mature during the inter-season to form the initial consumer and resource populations of season $(i+1)$. Up to some proportionnality constants accounting for the efficiency of the reproduction processes as well as overwintering mortality, \cite{Akhmetzhanov2010} obtained the following relation between the number of consumers of season $(i+1)$ and the initial number of resources of season $(i+1)$:
$$
c_{i+1}= \alpha J_i,\quad n_{i+1}(t=0)=\beta J_{n,i},
$$
with $J_i$ and $J_{n,i}$ defined in equations (\ref{e:value0}).

In the present mutant invasion model, the total consumer population is structured into $(1-\eps_i)c_i$ residents and $\eps_i c_i$ mutants that have different reproduction strategies. Taking into account this structure and assuming that mutants' progeny is also composed of mutants, we get the following inter-seasonal part for the mutant invasion model.
$$
c_{r_{i+1}} = \alpha \tilde U_r(c_i,\eps_i,n_i,T),\quad c_{m_{i+1}} = \alpha \tilde U_m(c_i,\eps_i,n_i,T),\quad n_{i+1}=\beta \tilde V(c_i,\eps_i,n_i,T)
$$
where the values $\tilde U_r$, $\tilde U_m$ and $\tilde V$ denote here the number of eggs/seeds produced by each (sub-)population:
$$
\tilde U_r = (1-\eps_i)c_i\int_0^T (1-u_r(t))p_r(t)\,\mathrm{d}t,\quad \tilde U_m = \eps_ic_i\int_0^T (1-u_m(t))p_m(t)\,\mathrm{d}t,
$$
$$
\tilde V = \int_0^T n(t)\,\mathrm dt\,.
$$

$\tilde U_r$, $\tilde U_m$ and $\tilde V$  can be computed from the solution of the optimal control problem (\ref{e:value3}) with the dynamics (\ref{e:dyn3}). Their values depend on the strategy chosen by the mutant and the resident, the length of the season $T$, the values $c_i$ and $\eps_i$ and initial conditions which are $p_r(0)=p_m(0)=0$ and $n(0)=n_i$. For a particular case $\eps=0$ the values $\tilde U_r$ and $\tilde U_m$ were derived analytically in subsection~\ref{s:valuefunctions}. In the following we investigate numerically this model on typical example; in particular we are interested in the long term fate of the resident and mutant consumer populations.

\begin{figure}[p]
$$\hspace{-.7cm}\includegraphics[scale=1.]{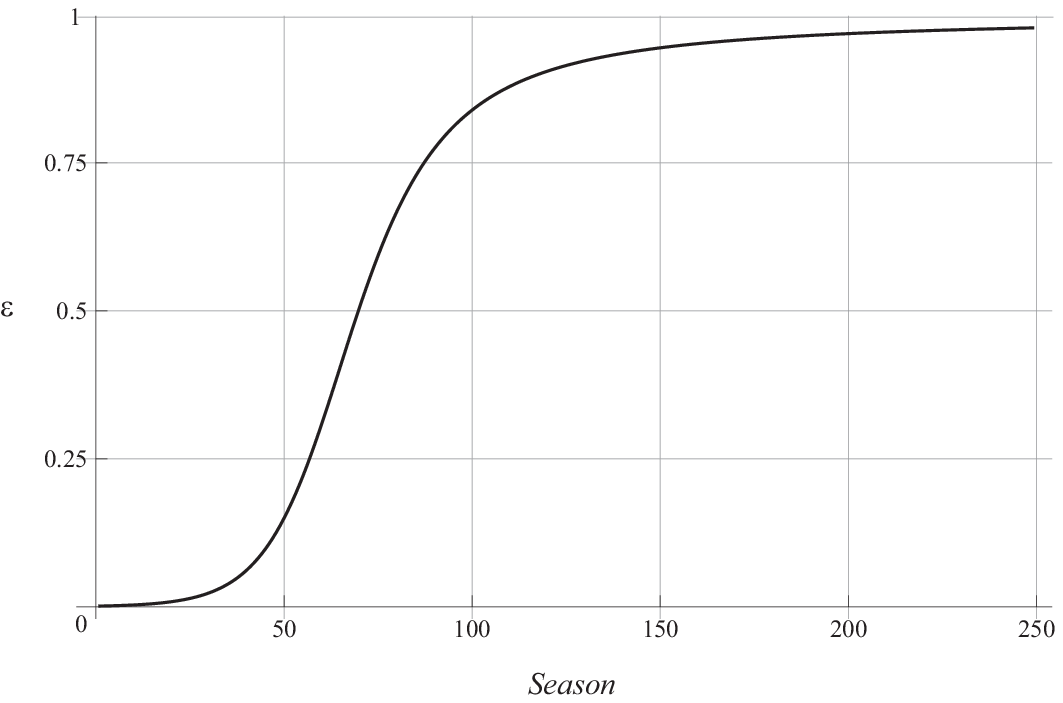}$$
\caption{Fraction of the mutants in the population of consumers}
\label{f:13}
$$\includegraphics[scale=1.]{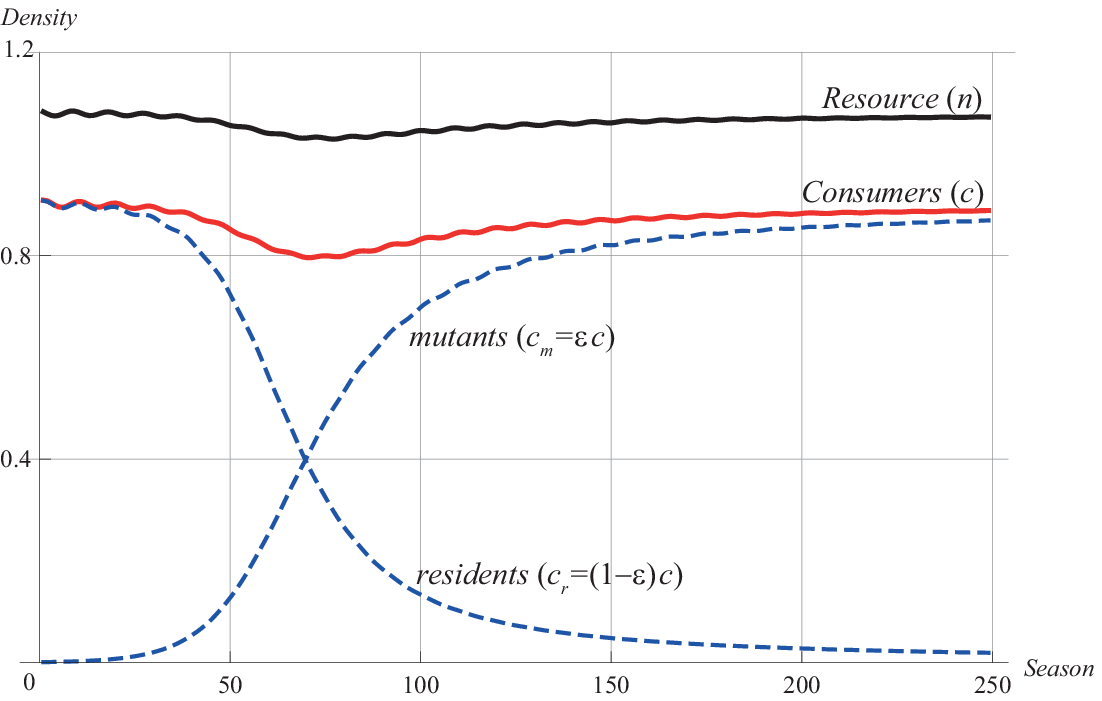}$$
\caption{Effect of the mutant invasion on the system}
\label{f:12}
\end{figure}

A previous investigation of the inter-seasonal model with collective optimal behavior of the consumers ({\it i.e.} there are no mutants such that $\eps=0$) has shown that the behavior of the system in long-term perspective could have very rich properties. Depending on the parameters of the problem, the value of $\beta$ and the length of season $T$, there could be an extinction of the resource or a or blowing up of the system (which leads to the suicide of the consumers). The system could also tend to some stable periodic behavior or to a globally asymptotic equilibrium. The last two cases illustrate a possible co-existence of the interacting species \cite{Akhmetzhanov2010}.

Here, we follow an adaptive dynamics like approach and consider that the resident consumer and the ressource population are at a (globally stable) equilibrium and investigate what happens when a small fraction of mutants appear in the resident consumer population. We actually assume that resident consumers are ``naive" in the sense that even if the mutant population becomes large through the season-to-season reproduction process, the resident consumers keep their collective optimal strategy and take mutants as cooperators, even if they do not cooperate.


We investigated numerically a case when $\alpha=2$, $\beta=0.5$ and $T=4$. The system is near the long-term  stable equilibrium point $c=0.9055$ and $n=1.0848$ as at the beginning of some season a mutant population of small size $c_m=0.001$ appears. The mutant population increases its frequency within the consumer population (see Fig.~\ref{f:13}) and modifies the dynamics of the system (Fig.~\ref{f:12}). The naive behavior of the consumers is detrimental to their progeny: along the seasons, mutant consumers progressively take the place of the collectively optimal residents and even replace them in the long run (Fig.~\ref{f:13}), making the mutation successful. We should however point out that the mutants' strategy as described in (\ref{e:value3}) is also a kind of ``collective" optimum: in some sense, it is assumed that the mutants cooperate with the other mutants. If the course of evolution drives the resident population to 0 and only mutants survive in the long run, this means that the former mutants become the new residents, with actually the exact same strategy as the one of the former residents they took the place from. Hence they are also prone to being invaded by non-cooperating mutants. The evolutionary dynamics of this naive resident-selfish mutant-resource appears thus to be a never-ending process: selfish mutants can invade and replace collective optimal consumers, but in the end transforms into collective optimal consumers as well, and a new selfish mutant invasion can start again. We are actually not in a ``Red Queen Dynamics" context since we focused on the evolution of one species only, and not co-evolution \cite{Valen1973}. Yet, what the Red Queen said to Alice seems to fit very well the situation we just described: ``here, you see, it takes all the running you can do to keep in the same place'' \cite{Carroll}.

\addcontentsline{toc}{section}{\textit{References}}
\bibliographystyle{unsrt}

\bibliography{./biblio_banff10}


\end{document}